\documentclass[12pt]{amsart}
\usepackage{amscd,amsmath,amsthm,amssymb,graphics}
\usepackage{amsfonts,amssymb,amscd,amsmath,enumerate,verbatim,calc,latexsym,pstricks,pst-plot,pst-3d,}
\usepackage{amsmath}
\usepackage{amssymb}
\usepackage{amsbsy}
\usepackage{graphicx}
\usepackage{amsfonts}
\usepackage{multicol}
\usepackage{amscd,amsmath,amsthm,amssymb}
\usepackage{pstricks,pst-plot,pst-3d}
\usepackage{lmodern,pst-node}
\usepackage{pstricks}

\usepackage{enumitem}

\makeatletter

\newcommand{\Rmnum}[1]{\expandafter\@slowromancap\romannumeral #1@}
\makeatother

\usepackage{pgf,tikz}
\usetikzlibrary{arrows}

\def\NZQ{\mathbb}               

\def\RR{{\NZQ R}}

\def\frk{\frak}               

\def\Phi{{\frk n}}
\def\Phi{{\frk N}}
\def\ub{{\mathbf u}}

\def\vb{{\mathbf v}}
\def\opn#1#2{\def#1{\operatorname{#2}}} 
\opn\chara{char} \opn\length{\ell} \opn\pd{pd} \opn\rk{rk}
\opn\projdim{proj\,dim} \opn\injdim{inj\,dim} \opn\rank{rank}
\opn\depth{depth} \opn\grade{grade} \opn\height{height}
\opn\embdim{emb\,dim} \opn\codim{codim}

\opn\Tr{Tr} \opn\bigrank{big\,rank}
\opn\superheight{superheight}\opn\lcm{lcm}
\opn\trdeg{tr\,deg}
\opn\reg{reg} \opn\lreg{lreg} \opn\ini{in} \opn\lpd{lpd}
\opn\size{size}\opn\bigsize{bigsize}
\opn\cosize{cosize}\opn\bigcosize{bigcosize}
\opn\sdepth{sdepth}\opn\sreg{sreg}
\opn\link{link}\opn\fdepth{fdepth}
\opn\index{index}
\opn\index{index}
\opn\indeg{indeg}
\opn\N{N}
\opn\SSC{SSC}
\opn\SC{SC}
\opn\conv{conv}
\opn\div{div} \opn\Div{Div} \opn\cl{cl} \opn\Cl{Cl}
\opn\Spec{Spec} \opn\Supp{Supp} \opn\supp{supp} \opn\Sing{Sing}
\opn\Ass{Ass} \opn\Min{Min}\opn\Mon{Mon} \opn\dstab{dstab} \opn\astab{astab}
\opn\Syz{Syz}
\opn\reg{reg}
\opn\Ann{Ann} \opn\Rad{Rad} \opn\Soc{Soc}
\opn\Im{Im} \opn\Ker{Ker} \opn\Coker{Coker} \opn\Am{Am}
\opn\Hom{Hom} \opn\Tor{Tor} \opn\Ext{Ext} \opn\End{End}
\opn\Aut{Aut} \opn\id{id}

\opn\nat{nat}
\opn\pff{pf}
\opn\Pf{Pf} \opn\GL{GL} \opn\SL{SL} \opn\mod{mod} \opn\ord{ord}
\opn\Gin{Gin} \opn\Hilb{Hilb}\opn\sort{sort}
\opn\initial{init}
\opn\ende{end}
\opn\height{height}
\opn\type{type}
\opn\aff{aff} \opn\con{conv} \opn\relint{relint} \opn\st{st}
\opn\lk{lk} \opn\cn{cn} \opn\core{core} \opn\vol{vol}
\opn\link{link} \opn\star{star}\opn\lex{lex}\opn\Mon{Mon}\opn\Min{Min}
\opn\gr{gr}

\def\pot#1#2{#1[\kern-0.28ex[#2]\kern-0.28ex]}

\opn\dirlim{\underrightarrow{\lim}}
\opn\inivlim{\underleftarrow{\lim}}
\def\Implies{\ifmmode\Longrightarrow \else
        \unskip${}\Longrightarrow{}$\ignorespaces\fi}
\def\implies{\ifmmode\Rightarrow \else
        \unskip${}\Rightarrow{}$\ignorespaces\fi}
\def\iff{\ifmmode\Longleftrightarrow \else
        \unskip${}\Longleftrightarrow{}$\ignorespaces\fi}

\let\:=\colon
\newtheorem{Theorem}{Theorem}[section]
 \newtheorem{Lemma}[Theorem]{Lemma}
 \newtheorem{Corollary}[Theorem]{Corollary}
 \newtheorem{Proposition}[Theorem]{Proposition}

 \newtheorem{Definition}[Theorem]{Definition}
 \newtheorem*{Definition*}{Definition}

 \newtheorem*{Conjecture*}{Conjecture}

 \newtheorem{Notation}[Theorem]{Notation}
 \newtheorem{Conventions}[Theorem]{Conventions}

\let\epsilon\varepsilon
\let\kappa=\varkappa
\def\qed{\ifhmode\textqed\fi
      \ifmmode\ifinner\quad\qedsymbol\else\dispqed\fi\fi}
\def\textqed{\unskip\nobreak\penalty50
       \hskip2em\hbox{}\nobreak\hfil\qedsymbol
       \parfillskip=0pt \finalhyphendemerits=0}
\def\dispqed{\rlap{\qquad\qedsymbol}}

\opn\dis{dis}
\def\pnt{{\raise0.5mm\hbox{\large\bf.}}}

\opn\Lex{Lex}

\begin{document}

 \title{Betti numbers of normal edge rings (\uppercase\expandafter {\romannumeral2})}

 \author{Zexin Wang}

\address{School  of Mathematical Sciences, Soochow University, 215006 Suzhou, P.R.China}
\email{zexinwang6@outlook.com}

 \author{Dancheng Lu*}

\footnote{* Corresponding author}

\address{School  of Mathematical Sciences, Soochow University, 215006 Suzhou, P.R.China}
\email{ludancheng@suda.edu.cn}

\begin{abstract}
We compute the Betti numbers of the edge rings of multi-path graphs using the \emph{induced-subgraph approach} introduced in \cite{WL1}. Here, a multi-path graph refers to a simple graph  composed of several paths that have the same starting point and the same ending point. Special cases include the graph $G_{r,d}$ introduced in \cite{GHK}, the graph $G_{r,s,d}$ introduced in \cite{NN}, and the graph $B_{\underline{\ell},h}$ introduced in \cite{LZ}. In particular, we show that all the multi-graded Betti numbers of a multi-path graph are the top multi-graded Betti numbers of some of its induced subgraphs.
\end{abstract}


 \maketitle 

\section{Introduction}

Let $G = (V(G), E(G))$ be a simple graph with vertex set $V(G) = \{x_1, \ldots, x_n\}$ and edge set $E(G) = \{e_1, \ldots, e_m\}$. Let $\mathbb{K}$ be a field. By definition, the \emph{edge ring} $\mathbb{K}[G]$ is the toric ring $\mathbb{K}[x_e \mid e \in E(G)] \subseteq \mathbb{K}[x_1, \ldots, x_n]$. Here, $x_e$ is the monomial  $\prod_{x_i \in e} x_i$ for all $e \in E(G)$. Let $\mathbb{K}[E(G)]$ or $\mathbb{K}[e_1, \ldots, e_m]$ be the polynomial ring in variables $e_1,\ldots,e_m$.  We can define a ring homomorphism
\[
\phi: \mathbb{K}[E(G)] \rightarrow \mathbb{K}[G]
\]
by letting $\phi(e_i) = x_{e_i}$ for all $i$. The \emph{toric ideal} of $\mathbb{K}[G]$ or $G$ is the ideal $I_G = \ker(\phi)$. It follows that $\mathbb{K}[G] \cong \mathbb{K}[E(G)] / I_G$. The study on  the connection between the algebraic invariants of $I_G$ (or $\mathbb{K}[G]$) and the combinatorial data associated with the graph $G$ is an active area of research in combinatorial commutative algebra. Relevant research can be found in \cite{BOV}, \cite{FHKT}, \cite{FKV}, \cite{GHK}, and \cite{NN}, among others.

  It is worth noting that $I_G$ admits two natural gradings. Since the polynomial ring  $\mathbb{K}[x_1,\ldots,x_n]$ has a natural $\mathbb{Z}_{\geq 0}^n$-grading,  its subring $\mathbb{K}[G]$  has  a  natural $\mathbb{Z}_{\geq 0}^n$-grading as well. It is  clear that there is  a natural $\mathbb{Z}_{\geq 0}^n$-grading on $\mathbb{K}[E(G)]$ such that $\phi$ is a multi-homogeneous  map, and thus  $I_G$ has a $\mathbb{Z}_{\geq 0}^n$-grading. Additionally, by assigning a degree of 1 to each edge $e$ in $E(G)$, both $\mathbb{K}[G]$ and $I_G$ also acquire a $\mathbb{Z}_{\geq 0}$-grading.

According to \cite[Chapter IV]{P}, both \(I_G\) and \(\mathbb{K}[G]\) admit both minimal \(\mathbb{Z}_{\geq0}^{n}\)-graded and \(\mathbb{Z}_{\geq0}\)-graded free resolutions over the  regular ring \(\mathbb{K}[E(G)]\). This fact allows us to establish the following relationships between their Betti numbers. First, we have \(\beta_{i,h}(I_G)=\beta_{i + 1,h}(\mathbb{K}[G])\) for all \(i\geq0\) and \(h\in\mathbb{Z}_{\geq0}^{n}\), and \(\beta_{i,j}(I_G)=\beta_{i + 1,j}(\mathbb{K}[G])\) for all \(i\geq0\) and \(j\in\mathbb{Z}_{\geq0}\).
Define \(|h|=\sum_{k = 1}^{n}h_k\) for \(h=(h_1,\ldots,h_n)\in\mathbb{Z}_{\geq0}^{n}\). Then the relationship between their \(\mathbb{Z}_{\geq0}^{n}\)-graded Betti numbers and their \(\mathbb{Z}_{\geq0}\)-graded Betti numbers is as follows:
\begin{equation}
\beta_{i,j}(\mathbb{K}[G]) = \sum_{\substack{h\in\mathbb{Z}_{\geq0}^{n} \\
|h| = 2j}} \beta_{i,h}(\mathbb{K}[G]).  \label{relation}
\end{equation}

Regarding the graded Betti numbers of edge rings, there are only several scattered  results.   In \cite[Theorem 5.1]{BOV}, the graded Betti numbers of $I_{\mathbf{K}_{2,d}}$ are obtained, where $\mathbf{K}_{2,d}$ is a complete bipartite graph with a bipartition $\{x_1,x_2\}$ and $\{y_1,\ldots,y_d\}$. In \cite{GHK}, Galetto et al. computed the graded Betti numbers of the toric ideal of $G_{r,d}$, where $G_{r,d}$ is the graph formed by adding a path of length $2r - 2$ between the two vertices $x_1,x_2$ of degree $d$ in $\mathbf{K}_{2,d}$.  Lately, Nandi and Nanduri \cite{NN} computed the graded Betti numbers of $G_{r,s,d}$, which is obtained  by adding a path of length $2s - 2$ between the two vertices of degree $d + 1$ in $G_{r,d}$. Let $\underline{\ell}=(\ell_1,\ldots,\ell_h)$, where each $\ell_i$ is a positive even number. In their work \cite{LZ}, Lu and Zhou introduced a class of graphs denoted by $B_{\underline{\ell},h}$.  These graphs consist of $h$ even paths of length $\ell_1,\ldots,\ell_h$ respectively, all sharing common starting and ending points. They computed the graded Betti numbers of the edge ring for these graphs $B_{\underline{\ell},h}$ under the condition that at most one even path has a length different from the others. Before the introduction of the induced-subgraph approach in \cite{WL1}, these classes of graphs are all instances for which the graded Betti numbers of their edge rings had been computed.  In view of the shared structural features of these graphs, we introduce the following definition.
\begin{Definition}\label{multi-path graph}\em
Let $m\geq 2$ and let $\underline{\ell} = (\ell_1, \ldots, \ell_m)$ be a given integral vector with $\ell_i\geq 1$ for all $i=1,\ldots,m$.   By definition, the graph $\mathbf{G}_{\underline{\ell}}$ is a simple graph consisting of $m$ paths of lengths $\ell_1, \ldots, \ell_m$ respectively, all of which  originate from a common vertex $v_1$ and converge at another shared vertex $v_2$, and  have no other common vertices except $v_1$ and $v_2$. We call such graphs \emph{multi-path} graphs (of type $\underline{\ell}$). Note that  there is at most one index $i$ such that $\ell_i=1$ since otherwise $G$ admits  multi-edges.
 \end{Definition}
It is clear that the graphs $\mathbf{K}_{2, d}$, $G_{s,d}$, $G_{r,s,d}$ and $B_{\underline{\ell},h}$ belong to the class of multi-path graphs. However, the latter is much larger since it is not bipartite in general.
In \cite{WL1}, we introduced the so-called induced-subgraph approach. With this approach, we evaluated the multi-graded Betti numbers of $I_G$ when $G$ is a two-ear graph or a compact graph. In this paper, we will apply the same approach to determine all the multi-graded Betti numbers, and consequently all the graded Betti numbers, of the edge ring associated with any multi-path graph.

To state  the induced-subgraph approach introduced in \cite{WL1},
 we first recall some necessary definitions and notations. By definition, a simple graph satisfies the \emph{odd-cycle} condition whenever any two minimal cycles $C_1$ and $C_2$ either share at least one vertex or are connected by an edge (which is) bridging a vertex from $C_1$ to a vertex from $C_2$. Recall that Simis et al. in \cite{SVV} and Ohsugi-Hibi in \cite{OH} independently established that the edge ring $\mathbb{K}[G]$ is a normal algebra precisely when the graph $G$ fulfills the odd-cycle condition. It is clear that every multi-path graph fulfills the odd-cycle condition.
 Assuming the normality of $\mathbb{K}[G]$, \cite[Theorem 6.31]{BG} reveals that the multi-graded canonical module of the multi-graded $\mathbb{K}[E(G)]$-algebra $\mathbb{K}[G]$, denoted as $\omega_{\mathbb{K}[G]}$, constitutes an ideal of $\mathbb{K}[G]$ with a generating set that possesses a geometric interpretation,   thereby simplifying its computation. Leveraging this understanding, we can determine the top multi-graded Betti numbers of $\mathbb{K}[G]$ by invoking Formula 6.6 from \cite{BG}. Note that the \emph{top multi-graded Betti numbers} of a finitely generated $\mathbb{Z}^{m}$-graded module $M$  refer to the numbers $\beta_{p,h}(M)$ for $h \in \mathbb{Z}^{m}$ if $M$ has a projective dimension $p$. We will also  utilize  the terminology the \emph{second top multi-graded Betti numbers} of $M$ to refer to the numbers $\beta_{p-1,h}(M)$ for $h \in \mathbb{Z}^{m}$.

We also need the following concept that was  given in   \cite{WL1}.
\begin{Definition} \label{BettiD}\em
Let $G$ be a connected graph that fulfills the odd-cycle condition. Denote by $p$ the projective dimension of  $\mathbb{K}[G]$. We call $G$ to be a \emph{top-Betti graph} provided that  if $\beta_{p,h}(\mathbb{K}[G])\neq 0$ for some $h\in \mathbb{Z}_{\geq 0}^{V(G)}$, then $\mathrm{supp}(h)=V(G)$. Here and after, for $h\in \mathbb{Z}_{\geq 0}^{V(G)}$, $\mathrm{supp}(h)$ represents  the set of $x_i\in V(G)$ with $h(x_i)\neq 0$.
\end{Definition}

The detailed steps of our approach is as follows. Let $G$ be a connected simple graph that satisfies the odd-cycle condition.
  \begin{enumerate}
 \item Initially, we compute  the projective dimension $p$ of $\mathbb{K}[G]$, and then identify a term order $<$ such that, for each $i$ with $1\leq i\leq p$, we can derive a suitable upper bound for the total Betti number $\beta_i(\mathbb{K}[E(G)]/\mathrm{in}_<(I_G))$.  In practice, the upper bounds are often achieved by demonstrating that $\mathrm{in}_<(I_G)$ has regular quotients.

  \item Secondly, we compute the top multi-graded Betti numbers $\beta_{p,h}(\mathbb{K}[G])$ of $\mathbb{K}[G]$ by determining the minimal generators of $\omega_{\mathbb{K}[G]}$.
  \item Next, for indices $i$ with $0<i < p$, we identify  all of the top-Betti  induced subgraphs $H$ of  $G$ such that  the projective dimension of $\mathbb{K}[H]$ is equal to $i$.
Subsequently, by leveraging  \cite[Formula 6.6]{BG}, we compute the top multi-graded Betti numbers of the edge rings associated with these induced subgraphs, and then sum them up.
  \item If the sum obtained in the third step matches the  upper bound for the total Betti number of $\mathbb{K}[E(G)]/\mathrm{in}_<(I_G)$ obtained the first step, then we have achieved our goal of computing the multi-graded Betti numbers $\beta_{i,h}(\mathbb{K}[G])$ for $h\in \mathbb{Z}_{\geq0}^{V(G)}$.  In this case,  $\beta_{i,h}(\mathbb{K}[G])$ is equal to the top Betti number $\beta_{i,h}(\mathbb{K}[H])$, where $H$ is the induced graph of $G$ on $\mathrm{supp}(h)$. Moreover, if $\beta_{i,h}(\mathbb{K}[G])\neq 0$, then the projective dimension of $\mathbb{K}[H]$ is equal to $i$. In addition, the  given upper bound for  $\beta_i(\mathbb{K}[E(G)]/\mathrm{in}_<(I_G))$ is actually the exact value of  $\beta_i(\mathbb{K}[E(G)]/\mathrm{in}_<(I_G))$.
  \item Otherwise, we  have to consider induced subgraphs  with edge rings of projective dimension $i+1$ and to compute the second top multi-graded Betti numbers of the edge rings of these graphs.
\end{enumerate}

In this paper, we will show  that by employing only the first four steps of this approach, one can compute all the multi-graded Betti numbers for the edge rings of multi-path graphs.  To introduce our main results, we need an additional notation.
\begin{Definition} \label{top-support} \em Let $H$ be  a simple graph satisfying the odd-cycle condition.
We define $\mathcal{N}_H$, which we call the \emph{top-support} of $H$, as the set of monomials $X^h$ for which $\beta_{p,h}(\mathbb{K}[H])$ is non-zero, where $p$ is the projective dimension of $\mathbb{K}[H]$.
\end{Definition}
 It is evident that the top multi-graded Betti numbers of $\mathbb{K}[H]$ are determined by the set $\mathcal{N}_H$.
The main result of this paper is as follows.

\begin{Theorem}[\emph{Theorem~\ref{main2}}]
Let $G$ be a multi-path graph $\mathbf{G}_{\underline{\ell}}$ with each $\ell_i > 1$. Let $\mathcal{N}_i(G)$ denote the union of all $\mathcal{N}_H$ such that $H$ is a top-Betti induced subgraph of $G$ satisfying $\mathrm{pdim}(\mathbb{K}[H]) = i$. Then there exists  a monomial order $<$ on $\mathbb{K}[E(G)]$ such that $\mathrm{in}_<(I_G)$ is square-free, and for all $i\geq 1$  and for each  monomial $\alpha\in \mathbb{K}[V(G)]$, we have
\[ \beta_{i,\alpha}(\mathbb{K}[G])=\beta_{i,\alpha}(\mathbb{K}[E(G)])/\mathrm{in}_{<}(I_G))
=
\begin{cases}
1, & \text{if } \alpha \in \mathcal{N}_i(G); \\
0, & \text{otherwise.}
\end{cases}
\]
\end{Theorem}
The paper \cite{CV} by Conca and Varbaro established that the extremal Betti numbers of an ideal are identical to those of its square-free initial ideal. Our current work, in view of the above results, conforms to this theoretical framework.

 We define the multi-path graph \(\mathbf{G}_{\underline{\ell}}\) to be  of \emph{even type} if all \(\ell_i\) are even, of \emph{odd type} if all \(\ell_i\) are odd, and of \textit{mixed type} otherwise. The graded Betti numbers for the edge rings associated with either an odd-type or an even-type multi-path graph $\mathbf{G}_{\underline{\ell}}$ are presented explicitly   in Corollary~\ref{oddgraded} and Corollary~\ref{4.6} respectively  in terms of numbers $\ell_1,\ldots,\ell_m$.

The structure of this paper is outlined as follows. First, in the second section, we introduce some fundamental concepts and results, which will serve as the foundation for our subsequent discussions. In Section 3, we classify multi-path graphs into odd type, even type, and mixed type, and show that the edge ring of a mixed type muti-path graph  can be expressed as  the tensor product of the edge rings of its even part and odd part. In the last section, we present the main result of this paper by computing the multi-graded Betti numbers of the edge rings associated with a multi-path graph.

\section{Preliminaries}

In this section, we provide a brief review of the notation and fundamental facts that will be utilized later on.

\subsection{Betti numbers}
Let $\mathbb{K}$ be a field, and let $R:=\mathbb{K}[X_1,\ldots,X_n]$ be the polynomial ring in variables $X_1,\ldots,X_n$ that are $\mathbb{Z}^m$-graded. Here, we do not require that $m=n$.   For a finitely generated $\mathbb{Z}^m$-graded  $R-$module $M$, there exists the minimal multi-graded free resolution of $M$ that has the form:
\begin{equation*}\label{free}0\rightarrow\underset{h\in\mathbb{Z}^m}{\bigoplus}R[-h]^{\beta_{p,h}(M)}\rightarrow \cdots \rightarrow \underset{h\in\mathbb{Z}^m}{\bigoplus}R[-h]^{\beta_{1,h}(M)}\rightarrow\underset{h\in\mathbb{Z}^m}{\bigoplus}R[-h]^{\beta_{0,h}(M)}\rightarrow M \rightarrow 0.\end{equation*}
Note that $R[-h]$ is the cyclic free $R$-module generated in degree $h$.
The {\it projective dimension} of $M$ are defined to be  $\mathrm{pdim}\,(M):={\mbox{max}}\,\{i\mid \beta_{i,\,h}(M)\neq 0 \mbox{ for some } h \}.$
The number $\beta_{i,h}(M):={\rm{dim}}_{\mathbb{K}}\mathrm{Tor}_i^R(M,\mathbb{K})_h$ is called the $(i,h)$-th {\it multi-graded Betti number} of $M$ and $\beta_{i}(M):=\sum_{h\in \mathbb{Z}^m}\beta_{i,h}(M)$ is called the $i$-th {\it total Betti number} of $M$. If $I$ is a homogeneous ideal of $R$, then the Betti number $\beta_{i,h}(R/I)$ equals $\beta_{i-1,h}(I)$ for all $i\geq 1$, making the computation of the Betti numbers for $I$ equivalent to that for its quotient ring.

Denote by $p$ the projective dimension of the module $M$. Subsequently, $\beta_p(M)$ is referred to as the {\it top total Betti number} of $M$, while $\beta_{p,h}(M)$ for $h\in \mathbb{Z}^m$ are designated as the {\it top multi-graded Betti numbers} of $M$. On the other hand, $\beta_{p-1}(M)$ is labeled as the {\it second top total Betti number} of $M$, and $\beta_{p-1,h}(M)$ for $h\in \mathbb{Z}^m$ are designated as the {\it second top multi-graded Betti numbers} of $M$.

If $R$ is a standard graded ring,  then $m=1$, and  the regularity of a finitely generated graded $R$-module $M$ is defined as
$\reg(M): = \max\left\{j - i \ |\ \beta_{i,j}(M) \neq 0\right\}.$

Let $I$ be a homogeneous ideal of $R$. If $I$ has a $\mathbb{Z}^m$-graded structure, then for any monomial order $<$ on $R$, the initial ideal $\mathrm{in}_<(I)$ is also a $\mathbb{Z}^m$-graded ideal of $R$.
It is well-known that the graded Betti numbers of $I$ are less than or equal to those of its initial ideal $\mathrm{in}_<(I)$ (see e.g. \cite[Corollary 3.3.3]{HH}). This result can be directly extended to the $\mathbb{Z}^m$-graded case of $I$, as shown in the following lemma. This lemma is necessary for us to determine the upper bound of the Betti numbers of edge rings.

\begin{Lemma}\label{total}Let $I$ be a homogeneous ideal of the polynomial ring $R=\mathbb{K}[X_1, \ldots, X_n]$ and also a $\mathbb{Z}^m$-graded ideal of $R$. For any monomial order $<$ on $R$ and all $i\geq 0$, $h\in \mathbb{Z}^m$, we have
$$\beta_{i,h}(R/I) \leq \beta_{i,h}(R/\mathrm{in}_<(I)).$$
In particular,
$$\beta_{i,j}(R/I) \leq \beta_{i,j}(R/\mathrm{in}_<(I)),\quad \beta_{i}(R/I) \leq \beta_{i}(R/\mathrm{in}_<(I)).$$
 \end{Lemma}
\subsection{Conventions and notions}  We  set conventions and  fix notions that play an important role in the rest  of this paper.

\begin{Conventions}\label{conventions}\em
Let $G$ be a simple graph, with its vertices labelled by letters or subscripted  letters, such as $v_1, v_2, u_{i,j}$, and so forth. We establish a correspondence between these labels and unit vectors in $\RR^{V(G)}$ as follows: $$v_1, v_2, u_{i,j}$$ correspond to $$\vb_1, \vb_2, \ub_{i,j}$$ respectively. With this notation, we can express $\RR^{V(G)}$ as
the set\begin{align*}
 \left\{ \sum_{i,j} a_{i,j}\ub_{i,j} + \sum_{i=1}^2 b\vb_i \mid a_{i,j}, b_i \in \RR\right\}.
\end{align*}
By way of example,  $\deg(v_1v_2)=\ub + \vb$, and $\deg(v_1^2u_{1,2}u_{1,3})=2\vb_1 + \ub_{1,2} + \ub_{1,3}$, and so on.

 We identify a monomial  with its multi-degree.  For example, we will  write $\beta_{i,uv}(\mathbb{K}[G])$ for $\beta_{i, \ub+\vb}(\mathbb{K}[G])$.
 Under this identification, we may look up $\mathbb{Z}_{\geq 0}^{V(H)}$ as a subset of $\mathbb{Z}_{\geq 0}^{V(G)}$ if $H$ is an induced subgraph of $G$.
\end{Conventions}

\begin{Notation}\label{notation} \em Let $G$ be a simple graph satisfying the odd-cycle condition. For every edge  $e=\{u,v\}$ of $G$, let $\mathbf{v}_e=\mathbf{u}+\mathbf{v}$. We  define $D_H$ to be the multi-degree $\sum_{e\in E(H)}\mathbf{v}_e\in \mathbb{Z}_{\geq 0}^{V(G)}$. Since we identify a monomial with its multi-degree, $D_G$ can alternatively be interpreted as the product of all edges of $H$. Here,  an edge is considered as the product of its two vertices. Furthermore, we use $\Theta_G$ to denote the product of all vertices with degree at least 2 in $G$.
\end{Notation}
It is straightforward to observe  that $\Theta_G^2$ divides $D_G$. Moreover, a simple graph $G$ satisfying the odd-cycle condition  is a top-Betti graph if and only if $\Theta_G$ divides every monomial in $\mathcal{N}_G$. Recall that  $\mathcal{N}_G$, defined in Definition~\ref{top-support}, denotes the top-support of $G$.

 According to Lemma~\ref{canonical}, we have
\begin{equation} \label{equation2} \mathcal{N}_G=\{D_G/\alpha\mid \alpha \mbox{ is a minimal generator of } \omega_{\mathbb{K}[G]}\}.\end{equation}
\subsection{Toric ideals of graphs}
Let $G$ be a simple graph, i.e., a finite graph without loops and multiple edges. Recall that a \emph{walk} of $G$ of length $q$ is a subgraph $W$ of $G$ such that $E(W)=\{\{i_0,i_1\},\{i_1,i_2\},\ldots, \{i_{q-1},i_q\}\}$, where $i_0,i_1,\ldots,i_q$ are vertices of $G$. An walk $W$ of $G$ is even if $q$ is even, and it is closed if $v_0=v_q$. The graph $G$ is \emph{connected} if, for any two vertices $i_j$ and $i_k$ of $G$, there is a walk between $i_j$ and $i_k$.  A \emph{cycle} with edge set $\{\{v_0,v_1\},\{v_1,v_2\},\{v_{q-1},v_q=v_0\}\}$ is a special even closed walk where $v_1,\ldots,v_q$ are pairwise distinct.  A cycle is called \emph{even} (resp. \emph{odd}) if $q$ is even (resp. odd).

The generators of the toric ideal of $I_G$ are binomials which is tightly related to even closed walks in $G$. Given an even closed walk $W$ of $G$ with
$$E(W)=\{\{v_0,v_1\},\{v_1,v_2\},\ldots,\{v_{2q-2},v_{2q-1}\},\{v_{2q-1},v_0\}\}$$
We associate $W$ with the binomial defined by
$$f_W:=\prod\limits_{j=1}^{q}e_{2j-1}-\prod\limits_{j=1}^{q}e_{2j},$$
where $e_j=\{v_{j-1},v_j\}$ for $1\leq j\leq 2q-1$ and $e_{2q}=\{v_{2q-1},v_0\}$. A binomial $f = u - v \in I_{G}$ is called a {\it primitive binomial} if there is no binomial $g = u' - v' \in I_{G}$ other than $f$ such that $u'|u$ and $v'|v$. An even closed walk $W$ of $G$ is a {\it primitive even closed walk} if its associated binomial $f_{W}$ is a primitive binomial in $I_G$.
 By e.g. \cite[Proposition 10.1.10]{V}, we have $\{f_W\:\; W \mbox{ is a primitive even closded walk}\}$ is a universal  Gr$\ddot{\mathrm{o}}$bner basis of $I_G$.
 In particular, it is a Gr$\ddot{\mathrm{o}}$bner basis of $I_G$ with respect to any monomial order. The set of primitive even walks of a graph $G$ was described in \cite{OH} explicitly.
\begin{Proposition} {\em (\cite[Lemma~3.2]{OH})} \label{primitive}
Let $G$ be a simple connected graph. Then a primitive even closed walk $\Gamma$ of $G$ is one of the following:
\begin{enumerate}
  \item $\Gamma$ is an even cycle of $G$;
  \item $\Gamma$=$(C_1,C_2)$, where $C_1$ and $C_2$ are odd cycles of $G$ having exactly one common vertex;
  \item $\Gamma$=$(C_1,\Gamma_1,C_2,\Gamma_2)$,  where $C_1$ and $C_2$ are odd cycles of $G$ having no common vertex and where $\Gamma_1$ and $\Gamma_2$ are walks of $G$ both of which combine a vertex $v_1$ of $C_1$ and a vertex $v_2$ of $C_2$.
\end{enumerate}
\end{Proposition}

 \subsection{Edge cones of simple graphs and canonical modules}
Let $G$ be a simple graph and denote  the  vertex set and edge set of $G$  by $V(G)$ and  $E(G)$ respectively.  A {\em matching} in $G$ refers to a subset $M \subset E(G)$ where each pair of distinct edges $e$ and $e'$ in $M$ satisfies $e \cap e' = \emptyset$. The {\em matching number} of $G$, denoted as $\mbox{mat}(G)$, represents the greatest possible cardinality among all matchings in $G$.

For a subset $W$ of $V(G)$, the {\it induced subgraph} $G_W$ is the graph with vertex set $W$ and for every pair $j,k\in W$, they are adjacent in $G_W$ if and only if they are adjacent in $G$. For a vertex $j\in V(G)$, let $G \setminus j$ be the induced subgraph $G_{V(G) \setminus \{j\}}$.
 A non-empty subset $T \subseteq V(G)$ is called an \textit{independent set} of $G$ if $\{j,k\} \not\in E(G)$ for any $j,k \in T$. For an independent set $T \subseteq V(G)$, we use $N_G(T)$  (or simply $N(T)$ for brevity) to  denote the neighborhood set of $T$, namely, $$N_G(T):=\{j\in V(G)\mid \{j,k\}\in E(G) \mbox{ for some  } k\in T\}.$$

 Let $G$ be a simple graph with vertex set $V(G)=\{1,\ldots,n\}$ and  edge set $E(G)$.
For any edge $f =\{i, j\} \in E(G)$, we define $\mathbf{v}_{f}= \mathbf{e}_{i}+ \mathbf{e}_{j}$,
where $\mathbf{e}_{i}$ represents the $i$th unit vector in $\mathbb{R}^{n}$. The edge cone of $G$, denoted as $\mathbb{R}_{+}(G)$, is the cone in $\mathbb{R}^{n}$ spanned
by $\{\mathbf{v}_{f} | f \in E(G)\}$. Specifically, it can be expressed as:
\[\mathbb{R}_{+}(G) := \left\{ \sum_{f \in E(G)}{a_{f}\mathbf{v}_{f}} \bigg| a_{f} \in \mathbb{R}_{\geq 0} \text{ for all } f \in E(G) \right\}.\]

Let us express the edge cone $\RR_{+}(G)$  in terms of linear inequalities, or equivalently, through its facets. Given a vertex $i$ of $G$, let $H_i, H_i^{+}$ and $H_i^{>}$ denote the following equalities or inequalities respectively:
\begin{align*} H_i: x_i=0,  \qquad H_i^+: x_i\geq0, \qquad
H_i^{>}:x_i>0.
\end{align*}

Given an independent set $T$ of $G$, let $H_T, H_T^{-}$ and $H_T^{<}$ denote the following equalities or inequalities respectively:
\begin{align*} H_T: \sum_{i\in T}x_i=\sum_{i\in N(T)}x_i,  \qquad H_T^-: \sum_{i\in T}x_i\leq\sum_{i\in N(T)}x_i, \qquad
H_T^{<}:\sum_{i\in T}x_i<\sum_{i\in N(T)}x_i.
\end{align*}
Note that $H_i\neq H_{\{i\}}$, since the latter represents the equality $x_i=\sum_{j\in N(i)}x_j$.  If $G$ is bipartite, the facets of $\RR_{+}(G)$ are given in \cite[Proposition 3.6 ]{VV}, that we record  in the following result for convenience.

\begin{Lemma}  \label{7} If $G$ is a bipartite connected graph with bipartition $(V_1, V_2)$,  then $F$ is a facet of $\RR_+(G)$ if and only if either
  \begin{enumerate}
      \item [$(a)$] $F=H_i\cap\RR_+(G)$ for some $i\in V(G)$ with $G\setminus i$ connected ;
      \vspace {1mm}
      \item [$(b)$] $F=H_{T}\cap\RR_+(G)$ for some $T\varsubsetneq V_1$, such that $G_{T\cup N(T)}$ and $G_{V(G)\setminus (T\cup N(T))}$ connected.
    \end{enumerate}
\end{Lemma}

Let $\mathcal{I}$ denote the set of all vertices satisfying condition (a) and $\mathcal{T}$ denote the set of all independent sets satisfying (b). Then, since $\mathrm{aff}(\mathbb{R}_+(G))=H_{V_1}=H_{V_2}$, we have
\begin{equation} \label{eq:ineq1}\tag{$\Delta$}
\mathbb{R}_+(G)=H_{V_1}\cap \bigcap_{i\in \mathcal{I}}H_{i}^{+} \cap \bigcap_{T\in \mathcal{T}}H_{T}^{-}
\end{equation}
and
 $$\mathrm{relint}(\mathbb{R}_+(G))=H_{V_1}\cap \bigcap_{i\in \mathcal{I}}H_{i}^>\cap \bigcap_{T\in \mathcal{T}}H_{T}^< .$$

It is  established by \cite[Corollary 2.3]{OH} that if $G$ satisfies the odd-cycle condition, then the edge ring $\mathbb{K}[G]$ is normal. Let  $\mathcal{S}(G)$ represent the semigroup composed of all $h\in \mathbb{Z}^n$ with $X^{h}\in \mathbb{K}[G].$
 According to \cite[Theorem 6.31]{BG}, the ideal of $\mathbb{K}[G]$ generated by all monomials $X^{h}$ with $h \in \text{relint}(\mathbb{R}{+}(G)) \cap \mathcal{S}(G)$ constitutes the $\mathbb{Z}^{V(G)}$-graded canonical module of the   $\mathbb{Z}^{V(G)}$-graded $^*$local algebra $\mathbb{K}[G]$, which we denote by $\omega_{\mathbb{K}[G]}$.   Then, in light of \cite[Formula 6.6]{BG}, we obtain the following result.

\begin{Lemma}\label{canonical} Let $G$ be a  simple graph  satisfying the odd-cycle condition and denote by $p$ the projective dimension of $\mathbb{K}[G]$.  Then $$\beta_{p-i,h}(\mathbb{K}[G])=\beta_{i,D_G-h}(\omega_{\mathbb{K}[G]})\qquad  \mbox{ for all }\quad 0\leq i\leq p\mbox{ and for all }h\in \mathbb{Z}_{\geq 0}^{V(G)}.$$
\end{Lemma}

The following  lemma, which is actually a special case of \cite[Corollary 2.5]{OHH},  serves as a cornerstone of the induced-subgraph approach.

\begin{Lemma}\label{start} Let $H$ be  an induced subgraph of a simple graph $G$. Then for $h\in \mathbb{Z}_{\geq0}^{V(G)}$ with $\mathrm{supp}(h)\subseteq V(H)$, we have $$\beta_{i,h}(\mathbb{K}[H])=\beta_{i,h}(\mathbb{K}[G]) \qquad \mbox{  for all }i\geq 0.$$
\end{Lemma}

Another well-known lemma that will be instrumental in our study concerns the projective dimension of the edge ring. For a rigorous proof of this lemma, we refer the readers to Lemma 2.9 of \cite{WL1}.
\begin{Lemma} \label{pd}
Let $G$ be a connected simple graph with vertex set $V(G)$ and edge set $E(G)$. If the toric ideal $I_G$ of $G$ has a square-free initial ideal $J$, then
\[
\mathrm{pdim}(\mathbb{K}[G]) =\mathrm{pdim}(\mathbb{K}[E(G)]/J) = \left\{
\begin{array}{ll}
|E(G)| - |V(G)| + 1, & \text{if $G$ is bipartite;} \\
|E(G)| - |V(G)|, & \text{otherwise.}
\end{array}
\right.
\]
\end{Lemma}

\section{Classifications of multi-graphs}

In this section, we classify multi-path graphs into odd type, even type, and mixed type. We observe that for a multi-path graph of mixed type, its edge ring is the tensor product of the edge rings of its two induced subgraphs,  which are of odd-type and even-type, respectively. In addition, all the top-Betti induced subgraphs of a multi-path graph are identified.

\begin{Definition} \label{multi-path}\em Let \(m\geq2\) be an integer, and \(\underline{\ell}=(\ell_1,\ldots,\ell_m)\) be a vector in \(\mathbb{Z}_{>0}^m\). We classify the multi-path graph \(\mathbf{G}_{\underline{\ell}}\) as follows: it is of \emph{even type} if all \(\ell_i\) are even, of \emph{odd type} if all \(\ell_i\) are odd, and of \textit{mixed type} if it does not fall into either the even-type or odd-type categories.

Assume that $G=\mathbf{G}_{\underline{\ell}}$ is of mixed type,  where $\ell_1,\ldots,\ell_s$ are even and $\ell_{s+1},\ldots,\ell_m$ are    odd. Define the {\it even part} of $G$ to be the induced subgraph of $G$, which is isomorphic to the  multi-path graph $\mathbf{G}_{(\ell_1,\ldots,\ell_s)}$. Similarly, define the {\it  odd part} of $G$ to be the induced subgraph of $G$, which is isomorphic to $\mathbf{G}_{(\ell_{s+1},\ldots,\ell_m)}$. The even part and the odd part  of $G$ are denoted by $G_e$ and $G_o$, respectively.

\end{Definition}

The exploration of the edge ring of a multi-path graph can be reduced to the investigation of the edge rings of multi-path graphs, which are of either even type or odd type, as  evidenced by the following lemma.

 \begin{Lemma} \label{reduction}Let $G=\mathbf{G}_{\underline{\ell}}$ be a multi-path graph of mixed type. Then   the edge ring of $G$ decomposes as the tensor product of the edge rings of $G_e$
  and $G_o$   over the base field $\mathbb{K}$, i.e.,$$\mathbb{K}[G]\cong \mathbb{K}[G_e]\bigotimes_{\mathbb{K}} \mathbb{K}[G_o].$$
  Moreover, for any monomial order $<_e$ on $\mathbb{K}[E(G_e)]$ and any monomial order $<_o$ on $\mathbb{K}[E(G_o)]$, there exists a monomial order $<$ on $\mathbb{K}[E(G)]$ such that $$\mathbb{K}[E(G)]/\mathrm{in}_<(I_G)\cong \left (\mathbb{K}[E(G_e)]/\mathrm{in}_{<_e}(I_{G_e})\right)\bigotimes_\mathbb{K}(\mathbb{K}[E(G_o)]/\mathrm{in}_{<_o}(I_{G_o})). $$
\end{Lemma}

\begin{proof} First, we note that $E(G)$ is the disjoint union of $E(G_e)$ and $E(G_o)$.   Next, since every pair of odd cycles of $G$ shares two common vertices,  by \cite[Lemma~3.2]{OH1} or Proposition~\ref{primitive},   primitive even closed walks of $G$  are exactly even cycles composed  of either two even paths or two odd paths of $G$. Consequently,  a primitive even closed walk of $G$  is contained in either $G_e$ or $G_o$. It follows from \cite[Proposition 10.1.10]{V} that a universal Gr\"obner basis of $I_G$ is the disjoint union of  a universal Gr\"obner basis of $I_{G_e}$ and  a universal Gr\"obner basis of $I_{G_o}$. These facts lead to the  desired decompositions.
\end{proof}

 According to Lemma~\ref{reduction}, the minimal multi-graded  free resolution of $\mathbb{K}[G]$ is the tensor production of the minimal multi-graded free resolutions  of $\mathbb{K}[G_o]$ and $\mathbb{K}[G_e]$.  As a direct consequence, we have for any $i\geq 0$ and any monomial $\alpha\in \mathbb{K}[G]$, the following formula holds: \begin{equation}\tag{*}\label{tensor}\beta_{i,\alpha}(\mathbb{K}[G])=\sum_{j+k=i}^{\beta\gamma=\alpha}\beta_{j,\beta}(\mathbb{K}[G_e])\cdot \beta_{k,\gamma}(\mathbb{K}[G_o])
\end{equation}
Similarly, if we replace $\mathbb{K}[G], \mathbb{K}[G_e], \mathbb{K}[G_o]$ by  $\mathbb{K}[E(G)]/\mathrm{in}_<(I_G), \mathbb{K}[E(G_e)]/\mathrm{in}_{<_e}(I_{G_e})$ and $ \mathbb{K}[E(G_o)]/\mathrm{in}_{<_o}(I_{G_o})$ respectively, the above formula remains hold.

\begin{Corollary} \label{reduction1} Let $G$ be a multi-path graph of mixed type. Then the following assertions hold.
\begin{enumerate}
\item [$\mathrm{(1)}$] $\mathrm{pdim} (\mathbb{K}[G])=\mathrm{pdim}(\mathbb{K}[G_e])+\mathrm{pdim}(\mathbb{K}[G_o])$;

\item [$\mathrm{(2)}$] $\mathcal{N}_G=\mathcal{N}_{G_e}\cdot \mathcal{N}_{G_o} \triangleq\{\alpha\beta\:\; \alpha\in \mathcal{N}_{G_e} \mbox{ and  } \mathcal{N}_{G_o}\}$;

\item  [$\mathrm{(3)}$] If  both $G_e$ and $G_o$ are top-Betti graph then so is $G$.
\end{enumerate}
\begin{proof} (1) and (2) follow from  formula~(\ref{tensor}). Note that for a simple graph $H$ that fulfills the odd-cycle condition, $H$ is a top-Betti graph if and only if for each $\alpha\in \mathcal{N}_H$, one has $\mathrm{supp}(\alpha)=V(H)$. Thus,   (3) follows from (2).  \end{proof}
\end{Corollary}

\begin{Proposition} \label{pdim} Let $G$ be a multi-path graph that consists of $m$ odd-length paths and $n$ even-length paths. Then
the projective dimension of $\mathbb{K}[G]$ equals  $m+n-2$ if  $m$ and $n$ are both $\geq 1$, otherwise it is $m+n-1$.

\begin{proof} Let $H$ be a multi-graph of even type or odd type and assume that $H$ consists of paths of lengths $\ell_1,\ldots,\ell_r$. Then $H$ is a bipartite graph with $|E(G)|=\ell_1+\cdots+\ell_r$ and $|V(G)|=(\ell_1-1)+\cdots+(\ell_r-1)+2$. It follows from Proposition~\ref{pd} that $\mathrm{pdim}(\mathbb{K}[H])=r-1$. Now, the result follows from (1) of Corollary~\ref{reduction1}.
\end{proof}
\end{Proposition}

We next discern all the top-Betti induced subgraphs of a multi-path graph.  Recall that in \cite[Proposition 3.1]{WL1} we showed that for a connected bipartite graph $G$, $G$ is a top-Betti graph if and only if for any $x\in V(G)$, one has $c_0(G\setminus x)< \deg(x)$. Moreover, it was shown  in \cite[Corollary 3.2]{WL1} that every vertex of a top-Betti graph has a degree at least 2. Here, $c_0(G)$ denotes the number of bipartite connected components of $G$.

\begin{Proposition}\label{top-betti}
Let $G$ be a multi-path graph and let $H$ be an induced subgraph of $G$. Then $H$ is a top-Betti graph if and only if the following conditions hold:
\begin{enumerate}
\item [$\mathrm{(1)}$] $H$ is also a multi-path graph;

\item [$\mathrm{(2)}$] If we let  $m_H$ and $n_H$ denote the number of even paths and odd paths in $H$ respectively, then both $m_H\neq 1$ and $n_H\neq 1$.
\end{enumerate}
\end{Proposition}
\begin{proof}  By \cite[Corollary 3.2]{WL1}, if $H$ is top-Betti, then $H$ does not contain vertices of degree one. This implies that $H$ is also a multi-path graph. Therefore, we may assume that $H$ is a multi-path graph. Let $H_e$ and $H_o$ be the odd part and the even part of $H$ respectively. Note that both $H_e$ and $H_o$ are bipartite.
Then, according to \cite[Proposition 3.1]{WL1}, $H_e$ is top-Betti if and only if $m_H\neq 1$, and $H_o$ is top-Betti if and only if $n_H\neq 1$. So, by Corollary~\ref{reduction1}.(3), if $m_H\neq 1$ and $n_H\neq 1$, then $H$ is top-Betti.
Conversely, if $m_H = 1$, then $I_{H_e}=0$, and thus $\mathcal{N}_{H_e}=\{0\}$. By Corollary~\ref{reduction1}.(2), we can conclude that $H$ is not top-Betti. The proof for the case when $n_H = 1$ is similar.
\end{proof}

\section{Multi-graded numbers of multi-path graphs}
In this section, we provide an explicit multi-graded Betti numbers formula for the edge rings of multi-path graphs by means of the induced-subgraph approach.
\subsection{Even type}
Let $G=\mathbf{G}_{2\underline{\ell}}$ be a multi-path graph of even type.  Here, $\underline{\ell}=(\ell_1,\ldots,\ell_m)\in \mathbb{Z}_{>0}^m$. We may assume the vertex set of $G$ is given by
$$V(G)=\{v_1,v_2\}\cup \{u_{i,j} \mid 1\leq i\leq m, 1\leq j\leq  2\ell_i-1 \}$$
while the edge set of $G$ is given by
\begin{align*}
E(G)=\{& \{v_1,u_{1,1}\},\ldots,\{v_1,u_{m,1}\}\}\cup   \{\{v_2,u_{1,2\ell_1-1}\},\ldots,\{v_2,u_{m,2\ell_m-1}\}\}
      \\ & \quad \cup   \{\{u_{i,j},u_{i,j+1}~|~1\leq i \leq m, 1\leq j \leq 2\ell_i-2\}\}.
\end{align*}
\noindent We label the edges of $G$ as follows: $e_{i,1}=\{v_1,u_{i,1}\},e_{i,j}=\{u_{i,j-1},u_{i,j}\},e_{i,2\ell_i}=\{v_2,u_{i,2\ell_i-1}\}$, for $i\in [m]$ and $j\in \{2,\ldots,2\ell_i-1\}.$

 It is evident that $G$ is a bipartite graph  with  bipartition $(V_{1}, V_{2})$, where
$$V_{1}=\bigcup\limits_{j=1}^{m}\bigcup\limits_{i=1}^{\ell_j}\{u_{j,2i-1}\}$$
and $$V_{2}=\{v_{1},v_{2}\}\cup \bigcup\limits_{j=1}^{m}\bigcup\limits_{i=1}^{\ell_j}\{u_{j,2i}\}.$$

Firstly, we establish an upper bound for the total Betti number of $\mathbb{K}[G]$, which is actually the precise value, as will be shown later.
For this, we need a notion that was originally introduced in \cite{M20} and subsequently explored in \cite{LZ}.
Recall that a monomial ideal $I$ in $R$ is said to possess {\it regular quotients} if there exists a linear ordering $\alpha_1 \prec \alpha_2 \prec \ldots \prec \alpha_r$ among the minimal generators $\{\alpha_1, \alpha_2, \ldots, \alpha_r\}$ of $I$ such that, for each $s=2,\ldots,r$, the quotient ideal $(\alpha_1, \ldots, \alpha_{s-1}):\alpha_s$ is generated by a regular sequence. Notably, a sequence of monomials $\beta_1, \ldots, \beta_r$ is regular precisely when the supports of distinct monomials are disjoint, i.e., $\mathrm{supp}(\beta_i) \cap \mathrm{supp}(\beta_j) = \emptyset$ for all $1 \leq s < t \leq r$. Here, the support of a monomial $\alpha$, denoted by $\mathrm{supp}(\alpha)$, refers to the set of variables that divide $\alpha$.

Assuming that the quotient ideal $(\alpha_1,\ldots,\alpha_{s-1}):\alpha_s$, denoted as $L_s$, is generated by a regular sequence of length $k_s$ for each $s$ ranging from $2$ to $r$, we can apply Koszul theory to deduce that the total Betti numbers $\beta_i(R/L_s)$ are given by $\beta_i(R/L_s) = \binom{k_s}{i}$ for all non-negative integers $i$. Let $I_s$ represent the monomial ideal generated by $\alpha_1,\ldots, \alpha_s$, and let $d_s$ represent the degree of $\alpha_s$ for each $s$ from $1$ to $r$.

Now, consider the following short exact sequences:
$$
0\rightarrow R/L_s[-d_s]\rightarrow R/I_{s-1}\rightarrow R/I_s\rightarrow 0, \quad \text{for all}\ s=2,\ldots,r.
$$
Utilizing the mapping cone construction (see e.g. \cite[section 27]{P} for the details), we can deduce that
$\beta_{i}(R/I_s)\leq \beta_{i-1}(R/L_s)+\beta_{i}(R/I_{s-1}).$
Repeatedly applying the above inequality, we arrive at
$$
\beta_i(R/I)\leq \beta_{i-1}(R/L_r)+\cdots+\beta_{i-1}(R/L_2)+\beta_{i}(R/(\alpha_1))
$$
and
$$
\beta_i(I)\leq \beta_{i}(R/L_r)+\cdots+\beta_{i}(R/L_2)+\beta_{i}(R).
$$
This implies that, for all non-negative integers $i$, we have the following inequality, which will be crucial in our subsequent discussions:
\begin{equation}\label{equation1}\beta_i(I)\leq\sum_{j=2}^r\binom{k_j}{i}+\binom{0}{i}.\end{equation}

\begin{Lemma}\label{even in} Let $G$  be the graph as described above. Then there exists a monomial order $<$ such that $\mathrm{in}_{<}(I_{G})$ is square-free and $\beta_{i}(\mathrm{in}_{<}(I_{G})) \leq (i+1)\binom{m}{i+2}$ for all $i \geq 0$. Thus, $$\beta_i(\mathbb{K}[G])\leq \beta_i(\mathbb{K}[E(G)]/\mathrm{in} _{<}(I_{G}))\leq i\binom {m}{i+1}.$$
\begin{proof}
According to Proposition 4.3 in \cite{LZ}, there exists a monomial order $<$ such that $\mathrm{in}_{<}(I_{G})$ is square-free with $\mathrm{in}_{<}(I_{G})=(u_1,\ldots, u_{\frac{m(m-1)}{2}})$, and the colon ideal $(u_1,\ldots, u_{j-1}):u_j$ is generated by a regular sequence of length $k_j$ for $j=2,\ldots,\frac{m(m-1)}{2}$, where
$$(k_2,\ldots,k_{\frac{m(m-1)}{2}})=\big(\underbrace{1,1}_{2},\underbrace{2,2,2}_{3},\ldots,\underbrace{m-2,\ldots,m-2}_{m-1}\big).$$
In view of Formula  (\ref{equation1}), we have,

\begin{align*}
\beta_{i}(\mathrm{in} _{<}(I_{G}))&\leq\sum_{k=1}^{m-2} (k+1)\binom{k}{i}+\binom{0}{i}
= (i+1)\sum_{k=1}^{m-2} \binom{k+1}{i+1}+\binom{0}{i}\\&
= (i+1)\sum_{k=0}^{m-2} \binom{k+1}{i+1}= (i+1)\binom{m}{i+2}.
\end{align*}
Here, the third-to-last equality follows from the identity $k\binom{n}{k} =
n\binom{n-1}{k-1}$, and the last equality follows from the identity
$\sum\limits_{j=0}^{n}\binom{j}{k}=\binom{n+1}{k+1}$ for $k\geq 0$.
From this bound and Lemma~\ref{total}, the second assertion is immediate.
\end{proof}
\end{Lemma}
The initial ideal $\mathrm{in}_{<}(I_{G})$ discussed in this subsection consistently refers to the initial ideal $\mathrm{in}_{<}(I_{G})$ studied in Lemma~\ref{even in}.

Next, we compute the minimal generators of the canonical module $\omega_{\mathbb{K}[G]}$. By Conventions~\ref{conventions}, we can express $\RR^{V(G)}$ as
the set\begin{align*}
 \left\{\sum\limits_{i=1}^{m}\sum\limits_{j=1}^{2\ell_{i}-1}{a_{i,j}\ub_{i,j}}+ b_{1}\vb_{1}+b_{2}\vb_{2} \mid a_{i,j}, b_{i}\in \RR \mbox{ for all } i,j \right\}.
\end{align*}
It is easy to see that  $|V(G)|=2\sum\limits_{i=1}^{m}\ell_{i}-m+2$. By \cite[Lemma 10.7.18]{V} which says $\RR_+(G)\cap \mathbb{Z}^{V(G)}=\mathcal{S}(G)$, we have $\mathrm{relint}(\RR_+(G))\cap \mathcal{S}(G) = \mathrm{relint}(\RR_+(G))\cap \mathbb{Z}^{V(G)}$. Subsequently, we construct $m-1$ integral vectors in $\RR^{V(G)}$ and then prove they are minimal in $\mathrm{relint}(\RR_+(G))\cap \mathbb{Z}^{V(G)}$.
Recall that an integral vector in $\mathrm{relint}(\RR_+(G))\cap \mathcal{S}(G)$ is  {\it minimal} if it cannot be written as the sum of a vector in $\mathrm{relint}(\RR_+(G))\cap \mathcal{S}(G)$ and a nonzero vector of  $\RR_+(G)\cap \mathcal{S}(G)$.

For $\ell=1,\ldots,m-1$, let $$\alpha_{\ell}:=\sum\limits_{i=1}^{m}\sum\limits_{j=1}^{2\ell_{i}-1}\ub_{i,j}+
\ell\vb_{1}+(m-\ell)\vb_{2}.$$

We now verify that $\alpha_{\ell} \in \mathrm{relint} (\RR_+(G))$ for every $1\leq \ell\leq m-1$.

First, by Lemma~\ref{7}, for a subset $T$ of $V_1$, we have $H_{T}\cap\RR_+(G)$ is the facet of $\RR_+(G)$ if and only if  $T$ is one of the following sets:
    \begin{enumerate}[label=(\roman*), leftmargin=*, nosep]
    \item $\displaystyle
        \bigcup_{j=1}^{s}\!\bigcup_{i=1}^{k_{p_j}} \{u_{p_{j},2i-1}\}, \quad
        \text{where } \{p_1,\ldots,p_s\} \subseteq [m], \; 1 \leq k_{p_j} \leq \ell_{p_j}-1$;
    \item $\displaystyle
        \bigcup_{j=1}^{s}\!\bigcup_{i=k_{p_j}}^{\ell_{p_j}} \{u_{p_{j},2i-1}\}, \quad
        \text{where } \{p_1,\ldots,p_s\} \subseteq [m], \; 2 \leq k_{p_j} \leq \ell_{p_j}$;
    \item $\displaystyle
        \bigcup_{\substack{j \neq p \\ i=1}}^{\ell_j} \{u_{j,2i-1}\}, \quad
        \text{where } p \in [m]$;
    \item $\displaystyle
        \left( \bigcup_{j \neq p} \bigcup_{i=1}^{\ell_j} \{u_{j,2i-1}\} \right) \cup
        \bigcup_{i=1}^{f_p} \{u_{p,2i-1}\}, \quad
        \text{where } p \in [m], \; 1 \leq f_p \leq \ell_p-1$;
    \item $\displaystyle
        \left( \bigcup_{j \neq p} \bigcup_{i=1}^{\ell_j} \{u_{j,2i-1}\} \right) \cup
        \bigcup_{i=k_p}^{\ell_p} \{u_{p,2i-1}\}, \quad
        \text{where } p \in [m], \; 2 \leq k_p \leq \ell_p$;
    \item $\displaystyle
        \left( \bigcup_{j \neq p} \bigcup_{i=1}^{\ell_j} \{u_{j,2i-1}\} \right) \cup
        \left( \bigcup_{i=1}^{f_p} \{u_{p,2i-1}\} \right) \cup
        \left( \bigcup_{i=k_p}^{\ell_p} \{u_{p,2i-1}\} \right) \quad
        \text{where } p \in [m], \; 1 \leq f_p < k_p-1 \leq \ell_p-1$;
    \item $\displaystyle
        \bigcup_{i=f_j}^{k_j} \{u_{j,2i-1}\}, \quad
        \text{where } j \in [m], \; 2 \leq f_j \leq k_j \leq \ell_j-1$.
\end{enumerate}
\allowdisplaybreaks
Hence, it follows from \eqref{eq:ineq1} that a vector of $\RR^{V(G)}$ of the form:
$\sum\limits_{i=1}^{m}\sum\limits_{j=1}^{2\ell_{i}-1}{a_{i,j}\ub_{i,j}}+b_{1}\vb_{1}+b_{2}\vb_{2}$ belongs to $ \RR_+(G)$ if and only if the following formulas are satisfied:
\begin{align*}
&(1) \quad \sum_{j = 1}^{m}\sum_{i = 1}^{\ell_j - 1}a_{j,2i}+b_{1}+b_{2} = \sum_{j = 1}^{m}\sum_{i = 1}^{\ell_j}a_{j,2i-1};\\
&(2) \quad a_{i,j} \geq 0, \quad \text{for any } 1\leq i\leq m \text{ and } 1\leq j\leq 2\ell_{i}-1;\\
&(3) \quad b_{1},b_{2} \geq 0;\\
&(4) \quad \sum_{j = 1}^{s}\sum_{i = 1}^{k_{p_{j}}}a_{p_{j},2i}+b_{1} \geq \sum_{j = 1}^{s}\sum_{i = 1}^{k_{p_{j}}}a_{p_{j},2i-1}, \quad \text{for any } \{p_1,\ldots,p_s\}\subseteq [m] \\ &\quad \text{ and for any } 1\leq k_{p_{j}}\leq \ell_{p_{j}}-1;\\
&(5) \quad \sum_{j = 1}^{s}\sum_{i = k_{p_{j}} - 1}^{\ell_{p_{j}} - 1}a_{p_{j},2i}+b_{2} \geq \sum_{j = 1}^{s}\sum_{i = k_{p_{j}}}^{\ell_{p_{j}}}a_{p_{j},2i-1}, \quad \text{for any } \{p_1,\ldots,p_s\}\subseteq [m] \\ &\quad \text{ and for any } 2\leq k_{p_{j}}\leq \ell_{p_{j}};\\
&(6) \quad \sum_{j\neq p}\sum_{i = 1}^{\ell_j - 1}a_{j,2i}+b_{1}+b_{2} \geq \sum_{j\neq p}\sum_{i = 1}^{\ell_j}a_{j,2i-1}, \quad \text{for any } p\in[m];\\
&(7) \quad \sum_{j\neq p}\sum_{i = 1}^{\ell_j - 1}a_{j,2i}+\sum_{i = 1}^{f_{p}}a_{p,2i}+b_{1}+b_{2} \geq \sum_{j\neq p}\sum_{i = 1}^{\ell_j}a_{j,2i-1}+\sum_{i = 1}^{f_{p}}a_{p,2i-1}, \\ &\quad\quad \text{for any } p\in[m] \text{ and } 1\leq f_p\leq\ell_p - 1;\\
&(8) \quad \sum_{j\neq p}\sum_{i = 1}^{\ell_j - 1}a_{j,2i}+\sum_{i = k_p}^{\ell_p}a_{p,2i - 2}+b_{1}+b_{2} \geq \sum_{j\neq p}\sum_{i = 1}^{\ell_j}a_{j,2i-1}+\sum_{i = k_p}^{\ell_{p}}a_{p,2i-1},\\&\quad \quad \text{for any } p\in[m] \text{ and } 2\leq k_p\leq\ell_p;\\
&(9) \quad \sum_{j\neq p}\sum_{i = 1}^{\ell_j - 1}a_{j,2i}+\sum_{i = 1}^{f_{p}}a_{p,2i}+\sum_{i = k_p}^{\ell_{p}}a_{p,2i - 2}+b_{1}+b_{2} \geq \sum_{j\neq p}\sum_{i = 1}^{\ell_j}a_{j,2i-1}+\sum_{i = 1}^{f_{p}}a_{p,2i-1} \\
&\quad\quad\quad +\sum_{i = k_p}^{\ell_{p}}a_{p,2i-1}, \qquad\text{for any } p\in[m] \text{ and } 1\leq f_p<k_p - 1\leq\ell_p - 1;\\
&(10) \quad \sum_{i = f_j - 1}^{k_j}a_{j,2i} \geq \sum_{i = f_j}^{k_j}a_{j,2i-1}, \quad \text{for any } j\in [m] \text{ and } 2\leq f_{j}\leq k_{j}\leq \ell_{j}-1.
\end{align*}

 For $1\leq \ell\leq m-1$, it is straightforward to check that the coordinates of $\alpha_{\ell}$ meet these inequalities in their strict forms. This implies that $\alpha_{\ell}$ belongs to $\mathrm{relint}(\RR_+(G))\cap \mathbb{Z}^{V(G)}$.

Next, we show that $\alpha_{\ell}$ is a minimal vector  in $\mathrm{relint}(\RR_+(G))\cap \mathbb{Z}^{V(G)}$.
Suppose on the contrary  that $\alpha_{\ell}=\alpha' + \alpha''$ for some $\alpha' \in \mathrm{relint}(\RR_+(G))\cap \mathbb{Z}^{V(G)}$ and $\alpha'' \in \RR_+(G)\cap \mathbb{Z}^{V(G)} \setminus \{{\bf 0}\}$.
Write
$$\alpha'=\sum\limits_{i=1}^{m}\sum\limits_{j=1}^{2\ell_{i}-1}{a'_{i,j}\ub_{i,j}}+b'_{1}\vb_{1}+b'_{2}\vb_{2},$$
$$\alpha''=\sum\limits_{i=1}^{m}\sum\limits_{j=1}^{2\ell_{i}-1}{a''_{i,j}\ub_{i,j}}+b''_{1}\vb_{1}+b''_{2}\vb_{2}.$$

In view of the inequalities (2) and (3), we see that $a'_{i,j}, b'_{i}\geq 1$ for all $i,j$. Hence, $a'_{i,j}= 1$ and $a''_{i,j}= 0$ for all possible $i,j$. Because of (6), we have $b'_{1} +b'_{2} \geq m$. Hence $b''_{1} +b''_{2} \leq0$.
From this together with (3) it follows that $b''_{1} =b''_{2} =0$.  Hence, $\alpha''=0.$ This is a contradiction, which shows that $\alpha_{\ell}$ is a minimal vector in $\mathrm{relint}(\RR_+(G))\cap \mathbb{Z}^{V(G)}$ for $\ell=1, \ldots, m-1$. We summarize  the above computation  in the following lemma.

Recall that $\Theta_G$ denotes the product of all vertices of degree at least 2 of $G$, see Notation~\ref{notation}.

\begin{Lemma}\label{generators} Let $G=\mathbf{G}_{2\underline{\ell}}$ be a multi-path graph of even type, where $\underline{\ell}=(\ell_1,\ldots,\ell_m)\in \mathbb{Z}_{>0}^m$. Then $\Theta_Gv_1^{\ell}v_2^{m-2-\ell}$ for $\ell=0,\ldots,m-2$ are minimal generators of $\omega_{\mathbb{K}[G]}$.
\end{Lemma}
\begin{proof} This is attributed to the existence of a bijection between the minimal vectors in $\mathrm{relint}(\RR_+(G))\cap \mathbb{Z}^{V(G)}$ and the minimal generators of $\omega_{\mathbb{K}[G]}$.
\end{proof}

Actually, $\Theta_Gv_1^{\ell}v_2^{m-2-\ell}$ for $\ell=0, \ldots, m-2$ are all minimal generators of $\omega_{\mathbb{K}[G]}$, as shown lately.

\begin{Lemma}\label{eventopsupp} Let $G$ be the graph described above.  Then $$\mathrm{pdim}(\mathbb{K}[G])=\mathrm{pdim}(\mathbb{K}[E(G)]/\mathrm{in}_<(I_G))=m-1.$$ Moreover, we have
  $$\mathcal{N}_G=\{\Theta_G v_{1}^{\ell}v_{2}^{m-2-\ell}\mid \ell=0,\ldots,m-2\}.$$
\end{Lemma}
\begin{proof}
The first assertion is directly obtained from Lemma~\ref{pd}.

Furthermore, since $\mathrm{pdim}(\mathbb{K}[G]) = m-1$, it follows from Lemma~\ref{even in} that the top Betti number of $\mathbb{K}[G]$ is bounded above by $m-1$. Consequently, the set $\mathcal{N}_{G}$ contains no more than $m-1$ elements. Conversely, Lemma~\ref{generators} asserts that $\mathcal{N}_{G}$ contains at least $m-1$ elements. Therefore, the second assertion is established by noting that $D_G = \Theta_G^2 (v_1v_2)^{m-2}$ in conjunction with Formula~(\ref{equation2}) in Subsection~2.3.
\end{proof}

\begin{Corollary} \label{evenreg}
For the graph $G$ defined earlier, the regularity of its edge ring $\mathbb{K}[G]$ is
\[
\mathrm{reg}(\mathbb{K}[G]) = \mathrm{mat}(G) - 1.
\]
 \end{Corollary}

\begin{proof}
By Corollary \ref{eventopsupp}, we have $\mathrm{pdim}(\mathbb{K}[G]) = m-1$ and all monomials in $\mathcal{N}_G$ have degree
\[
2\sum_{i=1}^m \ell_i - m + m = 2\sum_{i=1}^m \ell_i.
\]
 Furthermore, since $\mathbb{K}[G]$ is Cohen-Macaulay, its regularity satisfies
\[
\mathrm{reg}(\mathbb{K}[G]) = \max\{j \ |\ \beta_{m-1,j}(\mathbb{K}[G]) \neq 0\} - (m-1).
\]
Calculating this, we get
\[
\frac{1}{2} \cdot 2\sum_{i=1}^m \ell_i - (m-1) = \sum_{i=1}^m \ell_i - m + 1.
\]

One can directly verify that the edge set of $G$:
\[
\bigcup_{j=1}^{\ell_1} \left\{e_{1,2j-1}\right\} \cup \bigcup_{j=1}^{\ell_2} \left\{e_{2,2j}\right\} \cup \bigcup_{i=3}^{m} \bigcup_{j=1}^{\ell_i-1} \left\{e_{i,2j}\right\}
\]
forms a maximum matching, hence $\mathrm{mat}(G) = \sum_{i=1}^m \ell_i - m + 2$ as desired.
\end{proof}

\begin{Theorem}\label{even} Let $G$ be the graph defined above, and let $\mathrm{in}_<(I_G)$ denote the square-free initial ideal of the toric ideal $I_G$ studied in Lemma~\ref{even in}. Then, for any $0\leq i\leq m-1$, we have $$\beta_i(\mathbb{K}[G])=\beta_{i}(\mathbb{K}[E(G)]/\mathrm{in}_<(I_G)) .$$
 Moreover, if we let $\mathcal{N}_i(G)$ denote the set $\sqcup_{H} \mathcal{N}_H$, where the union is over all the multi-path induced subgraphs $H$ of $G$ that contain exactly $i+1$ even paths. Then for all $i\geq 1$ and  for any monomial $\alpha$ in $\mathbb{K}[V(G)]$, we have $$\beta_{i,\alpha}(\mathbb{K}[G])=\beta_{i,\alpha}(\mathbb{K}[E(G)]/\mathrm{in}_<(I_G))=\left\{
                                          \begin{array}{ll}
                                            1, & \hbox{$\alpha\in \mathcal{N}_i(G)$;} \\
                                            0, & \hbox{otherwise.}
                                          \end{array}
                                        \right.$$
 \end{Theorem}

\begin{proof} The case that $i=1$ is clear.  Fix $1\leq i\leq m-1$. By Propositions~\ref{pdim} and \ref{top-betti}, an induced subgraph $H$ of $G$ is a top-Betti graph whose edge rings has projective dimension $i$ if and only if $H$ is a multi-path graph   contain  exactly $i+1$ even paths of $G$. Since $G$ has $m$ even paths, there are $\binom{m}{i+1}$ such induced subgraphs.  Denoted  them  by $H_1,\ldots, H_r$, respectively.  Here, $r=\binom{m}{i+1}$. Note that $\beta_{i,\alpha}(\mathbb{K}[H_\ell])=1$ for any $\alpha\in \mathcal{N}_{H_\ell}$ by Lemma~\ref{eventopsupp}, it follows that $\beta_{i,\alpha}(\mathbb{K}[G])=1$ for any $\alpha\in \bigcup_{\ell=1}^r\mathcal{N}_{H_{\ell}}$ in view of Lemma~\ref{start}. Set $\mathcal{N}_i(G)=\bigcup_{\ell=1}^r\mathcal{N}_{H_{\ell}}$.

Note that $\mathrm{supp}(\alpha)=V(H_k)$ for any $\alpha\in \mathcal{N}_{H_k}$, we have $\mathcal{N}_{H_1},\ldots, \mathcal{N}_{H_r}$ are pairwise distinct. Hence,
$$|\mathcal{N}_i(G)|=\sum_{\ell=1}^r|\mathcal{N}_{H_\ell}|=ir=i\binom{m}{i+1}.$$ This implies $$i\binom{m}{i+1}=\sum_{\alpha\in \mathcal{N}_i(G)}\beta_{i,\alpha}(\mathbb{K}[G])\leq\beta_{i}(\mathbb{K}[G])\leq \beta_i(\mathbb{K}[E(G)]/\mathrm{in} _{<}(I_{G}))\leq i\binom{m}{i+1}.$$
The last inequality of the above inequality chain comes from Lemma \ref{even in}. Hence, $$\beta_{i}(\mathbb{K}[G])=\sum_{\alpha\in \mathcal{N}_i(G)}\beta_{i,\alpha}(\mathbb{K}[G]).$$ In particular, we have if $\beta_{i,\alpha}(\mathbb{K}[G])\neq 0$, then $\alpha\in \mathcal{N}_i(G)$ and so $\alpha\in \mathcal{N}_{H_k}$ for some $k\in [r]$. Now the conclusion follows by combining with Lemma \ref{total}.
\end{proof}

 In light of Theorem~\ref{even} together with formula~(\ref{relation}) in Introduction,  we are able to express the graded Betti numbers of $\mathbb{K}[G]$ in terms of the numbers $\ell_1,\ldots,\ell_m$.

\begin{Corollary} \label{4.6}  Let $G=\mathbf{G}_{2\underline{\ell}}$ be a multi-path graph of even type, where $\underline{\ell}=(\ell_1,\ldots,\ell_m)\in \mathbb{Z}_{>0}^m$.  For $1\leq i\leq m - 1$ and $j\geq 1$, we define $E_i^{j}$ to be the set
\[
\{(k_1,\ldots,k_{i + 1},r)\mid \ell_{k_1}+\cdots+\ell_{k_{i+1}}= j, 1\leq k_1<\cdots<k_{i+1}\leq m, 0\leq r\leq i-1\}.
\]
Then $\beta_{i,j}(\mathbb{K}[G])=\beta_{i,j}(\mathbb{K}[E(G)]/\mathrm{in}_<(I_G))$ is the cardinality of the set $E_i^{j}$ for all $i\geq 1$ and $j\geq 1$.
\end{Corollary}
We remark that the cardinality of $E_i^{j}$ is equal to $i$ times the cardinality of the set $\{(k_1,\ldots,k_{i + 1})\mid \ell_{k_1}+\cdots+\ell_{k_{i+1}}= j, 1\leq k_1<\cdots<k_{i+1}\leq m\}.$

\subsection{Odd type} Let \(G = \mathbf{G}_{2\underline{\ell}+1}\) be a multi-path graph of odd type, where \(\underline{\ell}=(\ell_1,\ldots,\ell_m)\in \mathbb{Z}_{\geq0}^m\). Note that there is at most one \(j\) such that \(\ell_j = 0\). If \(\ell_j = 0\) for some \(j\), then the Betti number of the edge ring of \(\mathbf{G}_{2\underline{\ell}+1}\) can be calculated using the Betti splitting,  for details see \cite[Corollary 3.12]{FHKT}. So it is reasonable to assume all \(\ell_i\geq 1\) from the beginning.  We remark that the results in this subsection could alternatively be derived using the induced-subgraph approach, in a manner analogous to the process for the even-type case. However, in the interest of simplifying the proofs, we opt not to directly employ the induced-subgraph approach. Instead, we will utilize the results from \cite{WL1} to obtain these findings.

We may assume that the vertex set of \(G\) is given by
\[V(G)=\{v_1,v_2\}\cup \{u_{i,j} \mid 1\leq i\leq m, 1\leq j\leq  2\ell_i \}\]
while the edge set of \(G\) is given by
\begin{align*}
E(G)=\{& \{v_1,u_{1,1}\},\ldots,\{v_1,u_{m,1}\}\}\cup   \{\{v_2,u_{1,2\ell_1}\},\ldots,\{v_2,u_{m,2\ell_m}\}\}\\
       & \quad \cup \{\{u_{i,j},u_{i,j + 1}\} \mid 1\leq i \leq m, 1\leq j \leq 2\ell_i-1\}.
\end{align*}

Recall that $I_G\subset \mathbb{K}[E(G)]$ admits a $\mathbb{Z}^{V(G)}$-grading.
    Note that the initial ideal $\mathrm{in}_{<}(I_G)$ of $I_G$  also admits a $\mathbb{Z}^{V(G)}$-grading with respect to any given order, since   every edge of $G$ has a $\mathbb{Z}^{V(G)}$-grading, i.e., the product of its two vertices. For every monomial $u\in \mathbb{K}[E(G)]$,  let $\mathbf{mdeg}(u)$ denote the  $\mathbb{Z}^{V(G)}$-degree of $u$, which can be looked upon as a monomial in $\mathbb{K}[V(G)]$.  We first show the following result.
\begin{Lemma} \label{new} Let $\alpha$ be a monomial in $\mathbb{K}[V(G)]$ such that $\beta_{i,\alpha}(I_G)\neq 0$ for some $i\geq 0$. Then there exists $\ell >0$ such that $\alpha=(v_1v_2)^\ell\beta$  and $v_i$ does not divides $\beta$ for $i=1,2$.
\end{Lemma}
\begin{proof}  In view of  \cite[Proposition 10.1.10]{V}, every primitive closed walk of $G$ is an even cycle consisting of two paths with common vertices $v_1$ and $v_2$. Moreover, the binomials associated with those cycle form a universal Gr\"onber base of $I_G$. Let $I$ be an initial ideal of $I_G$ with respect to some given monomial order. Then, for any minimal generator $u$ of $I$, we may write $u$ as $u=e_1^ue_2^uw_u$ such that $e_j^u$ is an edge of $G$ with $v_j$ dividing $\mathbf{mdeg}(e_j^u)$  for $j=1,2$. It is evident that, for any generators $u,v$ of $I$, we have either $\{e_1^u,e_2^u\}\cap \{e_1^v,e_2^v\}=\emptyset$ or $\{e_1^u,e_2^u\}=\{e_1^v,e_2^v\}$. By observing the Tayor resolution of $I$, we conclude that if $\beta_{i,\alpha}(I)\neq 0$, then  there exists $\ell >0$ such that $\alpha=(v_1v_2)^\ell\beta$  and $v_j$ does not divides $\beta$ for $j=1,2$. Since $\beta_{i,\alpha}(I_G)\leq \beta_{i,\alpha}(I)$(see Lemma \ref{total}), the result follows.
\end{proof}

 We require the concept of a compact graph, which is essential for our subsequent analysis. This concept was first introduced and categorized in \cite{WL}. A compact graph of type one, as stated in \cite[Theorem 3.11]{WL}, is a finite collection of odd cycles that share a common vertex.  Additionally, the multi-graded Betti numbers of the edge ring of a compact graph were determined in \cite{WL1}.

By identifying  $v_1$ with $v_2$ in the graph  $G$, we obtain a new  graph, which is actually a compact graph of type one. Denote  this new graph by $A$, and let $v$ be the new vertex resulting from the  identification  of  $v_1$ and $v_2$. Define a map $f: V(G)\rightarrow V(A)$ by $f(x)=x$ if $x\in V(G)\setminus \{x_1,x_2\}$ and $f(x)=v$ if $x\in \{v_1,v_2\}$. Then $\{x_1, x_2\}$ is an edge of $G$ if and only if $\{f(x_1),f(x_2)\}$ is an edge of $A$.

     The map:  $\{x_1, x_2\}\longmapsto \{f(x_1),f(x_2)\}$ establishes a bijection  between the edge sets $E(G)$ and $E(A)$. Moreover, through this bijection, the primitive even closed walks of $G$ correspond to the primitive even closed walks of $A$. It follows that the edge rings $\mathbb{K}[G]$ and $\mathbb{K}[A]$ possess identical defining ideals. Therefore, as graded rings, $\mathbb{K}[G]$ is isomorphic to $\mathbb{K}[A]$. However, as multi-graded rings, they are not isomorphic, as $\mathbb{K}[G]$ is $\mathbb{Z}^{V(G)}$-graded, while $\mathbb{K}[A]$ is $\mathbb{Z}^{V(A)}$-graded, and $|V(G)|=|V(A)|+1$. Thus, $K[G]$ and $K[A]$ admit the same minimal free resolution over $\mathbb{K}[E(G)]=\mathbb{K}[E(A)]$, but with different gradings. Let $u = x_1\cdots x_r$ be a monomial in $\mathbb{K}[V(G)]$, where each $x_i$ is a vertex of $G$, and $x_1,\ldots,x_r$ are allowed not to be distinct. We let $f(u)$ denote the monomial $f(x_1)\cdots f(x_r)\in\mathbb{K}[V(A)]$. Then, it is easy to see that $\beta_{i,u}(\mathbb{K}[A])=\sum_{f(v)=u}\beta_{i,v}(\mathbb{K}[G])$ for all $i\geq0$ and all monomials $u\in\mathbb{K}[V(A)]$. In light of this fact and Lemma~\ref{new}, we are able to adapt Lemma 3.2 and Theorem 4.11 on compact graphs of type one from \cite{WL1} to the case of multi-path graphs. Note that $\mathrm{mat}(G)=\mathrm{mat}(A)+1$.

\begin{Lemma} \label{5.3} If $G=\mathbf{G}_{2\underline{\ell}+1}$, then, $\mathcal{N}_G=\{\Theta_G v_1^{\ell}v_2^{\ell}\mid \ell=0,\ldots,m-2\}. $
\end{Lemma}

\begin{Theorem} \label{odd}
Let $G = \mathbf{G}_{2\underline{\ell} + 1}$ be a multi-path graph of odd type.
Then the following assertions hold:
\begin{enumerate}
\item[$(1)$]  $\mathrm{pdim}(\mathbb{K}[G]) = m-1$;
\item[$(2)$]  $\mathrm{reg}(\mathbb{K}[G]) = \mathrm{mat}(G)-1$;
\item[$(3)$] For each $i$ with $1 \leq i \leq m-1$, let $\mathcal{N}_i(G)$ denote the disjoint union of $\mathcal{N}_{H}$, where $H$ ranges through all multi-path induced subgraphs of $G$ consisting of $i+1$ paths.
Then there exists a monomial order $<$ on $\mathbb{K}[E(G)]$ such that $\mathrm{in}_<(I_G)$ is square-free, and for any monomial $\alpha\in \mathbb{K}[V(G)]$, we have
\[
\beta_{i,\alpha}(\mathbb{K}[G]) = \beta_{i,\alpha}(\mathbb{K}[E(G)]/\mathrm{in}_<(I_G)) =
\begin{cases}
1, & \text{if } \alpha \in \mathcal{N}_i(G), \\
0, & \text{otherwise}.
\end{cases}
\]
\end{enumerate}
\end{Theorem}

 We are able to express the graded Betti numbers of $\mathbb{K}[G]$ in terms of the numbers $\ell_1,\ldots,\ell_m$.

\begin{Corollary} \label{oddgraded}
Let $G = \mathbf{G}_{2\underline{\ell} + 1}$ be a multi-path graph of odd type, where $\underline{\ell}=(\ell_1,\ldots,\ell_m)$. For $1\leq i\leq m - 1$ and $j\geq 1$, we define $O_i^{j}$ to be the set
\[
\{(k_1,\ldots,k_{i + 1}, r)\mid \ell_{k_1}+\cdots+\ell_{k_{i+1}}+r = j, 1\leq k_1<\cdots<k_{i+1}\leq m, 1\leq r\leq i\}.
\]
Then $\beta_{i,j}(\mathbb{K}[G])=\beta_{i,j}(\mathbb{K}[E(G)]/\mathrm{in}_<(I_G))$ is the cardinality of the set $O_i^{j}$ for all $i\geq 1$ and $j\geq 1$.
\end{Corollary}

\subsection{General case}
In this subsection, we present the main result of this paper.
\begin{Lemma} \label{true-mixedsupp} Let $H$ be a multi-path graph  which consists of $m$ odd paths and $n$ even paths with $m,n\geq2$. Let $v_1$ and $v_2$ be vertices of $H$ with a degree $\geq 3$. Then $$\mathcal{N}_H= \{\Theta_Hv_1^{s+t+1}v_2^{n+s-t-1}\mid s=0,\ldots,m-2, t=0,\ldots,n-2\}.$$ In particular, $|\mathcal{N}_H|=(m-1)(n-1)$ and $\mathrm{supp}(\alpha)=V(H)$ for all $\alpha\in \mathcal{N}_H.$
 \end{Lemma}

Now, we are ready to present the main result of this paper.
\begin{Theorem} \label{main2} Let $G$ be a multi-path graph that consists of $m$ even-length paths and $n$ odd-length paths. Then
\begin{enumerate}

  \item [$(1)$]  For any integer $1\leq i\leq \mathrm{pdim}(\mathbb{K}[G])$, let $\mathcal{N}_i(G)$ denote the disjoint union of $\mathcal{N}_H$, where $H$ ranges over all top-Betti induced subgraphs of $G$ whose edge rings have a projective dimension equals to  $i$.
      Then there exists  a monomial order $<$ on $\mathbb{K}[E(G)]$ such that $\mathrm{in}_<(I_G)$ is square-free, and for all $i\geq 1$  and for each  monomial $\alpha\in \mathbb{K}[V(G)]$, we have
$$\beta_{i,\alpha}(\mathbb{K}[G])=\beta_{i,\alpha}(\mathbb{K}[E(G)]/\mathrm{in}_<(I_G))= \left\{
    \begin{array}{ll}
      1, & \hbox{$\alpha\in \mathcal{N}_i(G)$;} \\
      0, & \hbox{otherwise.}
    \end{array}
  \right.$$
    \item [$(2)$] Let $G_e$ and  $G_o$ denote the even part and odd part of $G$ respectively. Then,
  \[
\mathrm{reg}(\mathbb{K}[G]) = \left\{
\begin{array}{ll}
\mathrm{mat}(G) - 1 , & \text{if $G$ is even or odd or $m\geq 2$ and $n\geq 2$;} \\
\mathrm{mat}(G_e) - 1, & \text{if $n=1$ and $m\geq2$;}\\
\mathrm{mat}(G_o) - 1, & \text{if $m=1$ and $n\geq2$.}
\end{array}
\right.
\]
\end{enumerate}
\end{Theorem}
\begin{proof}
(1) From Proposition~\ref{primitive} and \cite[Proposition 10.1.10]{V}, we immediately see that when $n=1$, the multi-graded Betti numbers of $\mathbb{K}[G]$ coincide with those of $\mathbb{K}[G_e]$, and when $m=1$, they coincide with those of $\mathbb{K}[G_o]$. Hence, the cases in which either $n\in\{0,1\}$ or $m\in\{0,1\}$ are direct consequences of Theorem~\ref{even} and Theorem~\ref{odd}, respectively.
 Therefore, we may assume that $m \geq 2$ and $n \geq 2$. Fix $i\geq 1$. According to Lemma~\ref{reduction} together with  Theorem~\ref{even} and Theorem~\ref{odd}, we obtain
\begin{align*}
\beta_i(\mathbb{K}[G])=\beta_i(\mathbb{K}[E(G)]/\mathrm{in}_<(I_G))=\sum_{\substack{j+k=i \\ j \geq 0, k \geq 0}} \binom{m}{j+1} \binom{n}{k+1} j k
\end{align*}

By Lemma~\ref{total}, completing the proof reduces to showing that $|\mathcal{N}_i(G)| = \beta_i(\mathbb{K}[G])$. By Propositions~\ref{pdim} and \ref{top-betti}, there are three classes of top-Betti induced subgraphs whose edge rings have a projective dimension of $i$:
(i) the multi-path induced subgraphs of $G$ that contain $j \geq 2$ even paths and $k\geq 2$ odd paths with $j + k=i+2$;
(ii) the multi-path induced subgraphs of $G$ that contain $i + 1$ even paths and $0$ odd paths;
(iii) the multi-path induced subgraphs of $G$ that contain $0$ even paths and $i + 1$ odd paths.

Let $H_1,\ldots,H_r$ be the induced subgraphs of $G$ belonging to class (i), $H_{r+1},\ldots, H_s$ be the induced subgraphs belonging to class (ii), and $H_{s+1},\ldots, H_t$ be the induced subgraphs belonging to class (iii).

Fix integers $j$ and $k$ such that $j \geq 2$, $k \geq 2$, and $j + k=i + 2$. Then, there exist exactly $\binom{m}{j} \binom{n}{k}$ elements $p$ with  $1\leq p\leq r$ such that $H_p$ contains $j$ even paths and $k$ odd paths. Furthermore, for such $p \in [r]$, the cardinality of $\mathcal{N}_{H_p}$ is $(j - 1)(k - 1)$. This implies
\[
\left|\bigcup_{p = 1}^r \mathcal{N}_{H_p}\right|=\sum_{\substack{j + k=i + 2 \\ j \geq 2, k \geq 2}} \binom{m}{j} \binom{n}{k} (j - 1)(k - 1).
\]

Next, by Lemmas~\ref{5.3} and \ref{eventopsupp}, it is easy to see that $|\mathcal{N}_{H_p}|=i$ for $p=r + 1,\ldots,s,\ldots,t$. Moreover, we have $t - r=\binom{m}{i + 1}+\binom{n}{i + 1}$. It follows that
\[
\left|\bigcup_{p=r + 1}^t \mathcal{N}_{H_p}\right|=i\left(\binom{m}{i + 1}+\binom{n}{i + 1}\right).
\]

Therefore, we have
\[
|\mathcal{N}_i(G)| = \sum_{\substack{j + k=i + 2 \\ j \geq 2, k \geq 2}} \binom{m}{j} \binom{n}{k} (j - 1)(k - 1) + i\left(\binom{m}{i + 1}+\binom{n}{i + 1}\right)= \beta_i(\mathbb{K}[G]).
\]
This completes the proof.

(3) By the same reason as in (1), we may assume that $m\geq 2$ and $n\geq 2$. By Lemma \ref{true-mixedsupp}, the maximum degree of monomials in $\mathcal{N}_G$ is calculated as
\[
2\sum_{i=1}^m \ell_i + 2\sum_{i=m+1}^{m+n} \ell_i - n + n + 2(m-1) = 2\sum_{i=1}^{m+n} \ell_i + 2(m-1).
\]
Following analogous reasoning to the proof of Corollary \ref{evenreg}, we derive:
    \begin{align*}
       \mathrm{reg}(\mathbb{K}[G]) = & \frac{1}{2} \cdot \left(2\sum_{i=1}^{m+n} \ell_i + 2(m-1)\right) - \mathrm{pdim}(\mathbb{K}[G]) \\
                      =&\sum_{i=1}^{m+n} \ell_i - n + 1.
    \end{align*}
One can directly verify that the subset of $E(G)$:
\[
\bigcup_{i=1}^{m} \bigcup_{j=1}^{\ell_i} \left\{e_{i,2j}\right\} \cup \bigcup_{j=1}^{\ell_{m+1}} \left\{e_{m+1,2j-1}\right\} \cup \bigcup_{j=1}^{\ell_{m+2}} \left\{e_{m+2,2j}\right\} \cup \bigcup_{i=m+3}^{m+n} \bigcup_{j=1}^{\ell_i-1} \left\{e_{i,2j}\right\}
\]
forms a maximum matching, hence $\mathrm{mat}(G) = \sum_{i=1}^{m+n} \ell_i - n + 2$, as desired.
\end{proof}

\vspace{4mm}

{\bf \noindent Acknowledgment:} This project is supported by NSFC (No. 11971338).  The authors are grateful to the software systems \cite{C}
for providing us with a large number of examples that help us develop ideas and test  the results.

\end{document}